\documentclass[11pt,a4paper]{amsart}
\usepackage[all]{xy}
\usepackage{amssymb}

\numberwithin{equation}{section}

\newcommand{\ie}{{\em i.e.}\ }

\newtheorem{theorem}{Theorem}

\newtheorem{lemma}{Lemma}
\newtheorem{proposition}{Proposition}
\newtheorem{corollary}{Corollary}

\newcommand{\opname}[1]{\operatorname{\mathsf{#1}}}

\renewcommand{\mod}{\opname{mod}\nolimits}

\newcommand{\dimv}{\underline{\dim}\,}
\newcommand{\ev}{\underline{e}\,}\newcommand{\mv}{\underline{m}\,}
\newcommand{\fv}{\underline{f}\,}\newcommand{\gv}{\underline{g}\,}

\renewcommand{\ker}{\opname{ker}\nolimits}

\newcommand{\Z}{\mathbb{Z}}
\newcommand{\N}{\mathbb{N}}
\newcommand{\Q}{\mathbb{Q}}
\newcommand{\C}{\mathbb{C}}
\newcommand{\F}{\mathbb{F}}
\newcommand{\fq}{{\mathbb F}_q}

\newcommand{\Hom}{\opname{Hom}}

\newcommand{\gr}{\opname{Gr}}

\newcommand{\Ext}{\opname{Ext}}
\newcommand{\Aut}{\opname{Aut}}

\newcommand{\GL}{\opname{GL}}

\newcommand{\ca}{{\mathcal A}}

\newcommand{\cF}{{\mathcal F}}

\newcommand{\ch}{{\mathcal H}}

\newcommand{\co}{{\mathcal O}}

\newcommand{\hmod}{H\opname{-mod}}

\begin{document}

\title{On the quiver Grassmannian in the acyclic case\\}
\author{Philippe Caldero}
\address{Institut Camille Jordan, Universit\'e Claude Bernard Lyon 1,
69622 Villeurbanne Cedex, France}
\email{caldero@math.univ-lyon1.fr}
\author{Markus Reineke}
\address{Fachbereich C\\
   Bergische Universit\"at Wuppertal \\
   D - 42097 Wuppertal\\
   Germany}
\email{reineke@math.uni-wuppertal.de}

\begin{abstract}
Let $A$ be the path algebra of a quiver $Q$ with no oriented cycle. We study geometric properties of the Grassmannians of submodules of a given $A$-module $M$. In particular, we obtain some sufficient conditions for smoothness, polynomial cardinality and we give different approaches to Euler characteristics. Our main result is the positivity of Euler characteristics when $M$ is an exceptional module. This solves a conjecture of Fomin and Zelevinsky for acyclic cluster algebras.

\end{abstract}

\maketitle
\setcounter{section}{-1}
\section{Introduction}
Let $M$ be a finite dimensional space on a field $k$. The Grassmannian $\gr_e(M,k)$ of
$M$ is the set of subspaces of dimension $e$. It is well known that $\gr_e(M,k)$ is
an algebraic variety with nice properties. For instance, the linear group
$\GL_e(M,k)$ acts transitively on $\gr_e(M,k)$ with parabolic stabilizer, hence the
variety  $\gr_e(M,k)$ is smooth and projective. 

Suppose now that $M$ has a structure of $A$-module, where $A$ is a finitely
generated $k$-algebra. It is natural to define the Grassmannian $\gr_e(M,A)$ of
$A$-submodules of $M$ of given dimension $e$. It is a closed subvariety of
$\gr_e(M,k)$, hence it is a projective variety. But in general this variety is not
smooth.

Let $Q$ be a finite quiver with no oriented cycle and let $A:=kQ$ be the path
algebra of $Q$. Then, $A$ is a finite dimensional algebra with minimal idempotents
$e_i$, $i\in Q_0$, where $Q_0$ is the set of vertices of $Q$. The $A$-module $M$ has
a dimension vector $\underline{\dim}(M):=(\dim(e_iM))_i\in\N^{Q_0}$ which
refines the notion of dimension. The Grassmannian  $\gr_e(M,A)$ splits into a
disjoint union of closed subvarieties $\gr_{\ev}(M,A)$ of submodules of $M$ with
given dimension vector $\ev$ such that $\sum_i\ev_i=e$. Hence it is natural to study
the geometric properties  of the so-called quiver Grassmannian
$\gr_{\ev}(M):=\gr_{\ev}(M,A)$.

The quiver Grassmannian can be found in \cite{crawley,schofield} in connection with the dimension of spaces of morphims of $kQ$-modules.  It plays an important role in Fomin-Zelevinsky's  theory of cluster algebras, \cite{fominzelevinsky1}. Indeed, as proved in \cite{caldchap} and \cite{caldkell2}, Euler characteristics $\chi_c(\gr_{\ev}(M))$ of quiver Grassmannians define sets of generators of cluster algebras in the acyclic type. 

To be more precise, consider the ring of Laurent polynomial $\Z[x_i^{\pm},\,i\in Q_0]$ and for any (f.d.) $kQ$-module $M$ with dimension vector $\mv=(m_i)$, set:
$$X_M:=\prod_{i\in Q_0}x_i^{-m_i}\sum_{\ev}\chi_c(\gr_{\ev}(M))\prod_{h\in Q_1}x_{s(h)}^{m_{t(h)}-e_{t(h)}}x_{t(h)}^{e_{s(h)}},$$
where $Q_1$ is the set of arrows of $Q$, and $s(h)$, $t(h)$ are respectively the source and the target of the arrow $h$. Then, it is known that a set of generators of the cluster algebra  is given by $\{X_M\}$ where $M$ runs over the set of {\it exceptional modules}, \ie modules with no self-extension. These generators are called cluster variables. Moreover, the set $\{X_M\}$, where $M$ runs over the set of f.d. modules is (at least conjecturally) a $\Z$-base of the cluster algebra which is an analogue of Lusztig's
 dual semicanonical base, \cite{caldzel}.

 The aim of the paper is to study  geometric properties of quiver Grassmannians, and in particular to understand their Euler characteristics. 
 
 In the first section with give an inductive formula to compute $\chi_c(\gr_{\ev}(M))$ when $M$ is a preprojective (or postinjective) module. This uses the Auslander-Reiten algorithm of construction of the preprojective representations and in particular the notion of almost split sequences. Unfortunately, it is not clear from the inductive formula that the characteristics $\chi_c(\gr_{\ev}(M))$ are positive.

 The second section is devoted to positivity. On the one hand, a well-known result of J. A. Green and C. Ringel (see \cite{Green,Ringel}) asserts that quantum groups can be realized as deformed Hall algebras of the category $\mod kQ$ and so $\chi_c(\gr_{\ev}(M))$ can be realized as coefficients of some convolution product in a PBW-basis. On the other hand, we know that Lusztig's canonical basis of quantum groups has nice positivity properties. We obtain two results from these remarks: 
 
 First,  proposition \ref{hall-positive} asserts that if the quiver Grassmannian $\gr_{\ev}(M)\mid_{\fq}$ is counted by a polynomial of $q$, then its characteristics is positive.  
 
 Our second positivity result, theorem \ref{except},  solves a conjecture of Fomin and Zelevinsky for acyclic cluster algebras. It states that  $\chi_c(\gr_{\ev}(M))$ is positive if $M$ is an exceptional module.
 
 Recall that the positivity in the Dynkin case was proved in \cite{caldkell}. With the help of proposition \ref{hall-positive}, we generalize this result for extended Dynkin quivers in section 3: the Euler characteristics is positive for any representation $M$ on an extended Dynkin quiver. Actually, what we need  is a polynomiality property of the quiver Grassmannians $\gr_{\ev}(M)$ and we prove that we can deduce it from a polynomial property of   varieties of representations  which are easier to handle, namely the varieties of representations $X$ such that $\dim\Hom(X,M)$ is equal to a fixed integer $a$.
 
 In section 4, we give some results on smoothness. We realize the quiver Grassmannian as a geometric quotient of an affine variety and this enables us to obtain a description of the tangent space of $\gr_{\ev}(M)$ at a point  $U$; we deduce the following
 $$\langle\ev,\mv-\ev\rangle\leq\dim T_U(\gr_{\ev}(M))\leq \langle\ev,\mv-\ev\rangle+\dim\Ext^1(M,M),$$
 where $\langle\_,\_\rangle$ is the Euler form on $\Z^{Q_0}$ identified with the Grothendieck group of $\mod kQ$, see \ref{tangsp}. This inequality implies clearly that the quiver Grassmannian is smooth if the module $M$ is exceptional. This also provides a dimension estimate of the variety.
 
 We introduce in section 5 some material on cluster algebras in order to explain the main application of the positivity theorem. As explained above, cluster variables for acyclic quivers are positive, which is a particular case of a conjecture of Fomin and Zelevinsky.
 
 We end with an open problem tightly linked with the Euler characteristics and its positivity. Can we find a cellular decomposition of the quiver Grassmannian? Solving this problem can lead to some combinatorics for the calculus of $\chi_c(\gr_{\ev}(M))$.

 {\bf Acknowledgements:} The first author wants to thank Gus Lehrer and Andrei Zelevinsky for stimulating conversations.

\section{Quiver Grassmannians.}
\begin{subsection}{}
Let $k$ be a field. We fix a  finite quiver $Q$ such that the path algebra $H=kQ$ is finite dimensional. Since this property is equivalent to the fact that $Q$ has no oriented cycle, such a quiver will be called acyclic.

Let $Q_0$ be the set of vertices and
$Q_1$ the set arrows of $Q$. Set $n=\#Q_0$. For any arrow $h$ of $Q_1$, we denote by $s(h)$ its source and $t(h)$ its target. We fix an ordering of $Q_0=[1,n]$ such that if
there exists a (non trivial) path from $j$ to $i$, then $j<i$. Since $Q$ has no
oriented cycles, it is possible to construct such an ordering.

A representation $(V,x)$ of $Q$ over $k$ is a $Q_0$-graded
$k$-vector space $V=\oplus_{i\in Q_0} V_i$ together with an element
$x=(x_h)_{h\in Q_1}$ in $E_V:=\prod_{h\in Q_1}\Hom(V_{s(h)},V_{t(h)})$. 
A morphism between two representations $(V,x)$ and $(V',x')$ is a $Q_0$-uple of ($k$-spaces) morphisms in $\prod_{i\in Q_0}\Hom_k(V_i,V_i')$ compatible with $x$ and $x'$, {\it i.e.} inducing commuting square diagram.

We denote by Rep$_kQ$ the category of (finite dimensional) representations of $Q$ over $k$. For each vertex $i$, let $S_i$ be the associated simple representation. The Grothendieck group of the category Rep$_kQ$ can be identified with $\Z^{Q_0}$ via the dimension vector map $\dimv$:
$$\dimv(M)=(\dim(M_i))_{i\in Q_0}.$$

There is a well known equivalence between Rep$_kQ$ and the category $\mod(kQ)$ of finite dimensional $kQ$-modules; hence a representation of $Q$ will be naturally considered as a $kQ$-module.

The group $G_V:=\prod_{i\in Q_0} \GL(V_i)$ acts on $E_V$ by
$(g_i).(x_h)=(g_{t(h)}x_hg_{s(h)}^{-1})$.
Clearly, the isoclasses of
a finite-dimensional $H$-module $M$ are naturally identified with the points of a
$G_V$-orbit ${\mathcal O}_M$ of representations of $Q$.

For any finite dimensional $H$-module $M$, and any $\ev$ in  $\Z^{Q_0}$, we denote by
$\gr_{\ev}(M)$ the Grassmannian of submodules of $M$ with dimension
vector $\ev$:
\[\gr_{\ev}(M)=\{N,\,N\in\hmod,\,N\subset M,\, \dimv(N)=\ev\}.\]
It is a closed subvariety (hence projective) of the classical Grassmannian of the
vector space $M$. 

We know  that this variety
is obtained by base change from a variety defined over $\Z$. So the field $k$ will sometimes be omitted. Let $\chi_c$ be the Euler-Poincar\'e
characteristic of $l$-adic cohomology with proper support defined
by
$$\chi_c(\gr_{\ev}(M)) = \sum_{i=0}^\infty (-1)^i \dim H^i_c(\gr_{\ev}(M), \overline\Q_l).$$
\end{subsection}
\begin{subsection}{}
We give here some basic properties for the computation of the characteristics of quiver Grassmannians.
This first proposition asserts that we can restrict ourselves to the indecomposable case.
\begin{proposition}\label{split}\cite[Proposition 3.6]{caldchap}
Let $M$ and $N$ be two $Q$-representations and let ${\gv}$ be a dimension vector, then
$$\chi_c(\gr_{\gv}(M\oplus N))=\sum_{\ev+\fv=\gv}\chi_c(\gr_{\ev}(M))\chi_c(\gr_{\fv}(N)).$$
\end{proposition}
\begin{proof}
It is sufficient to remark, see \cite[Equation (26)]{caldchap}, that the map
$$\pi\,:\,\gr_{\gv}(M\oplus N)\rightarrow\coprod_{\ev+\fv=\gv}\gr_{\ev}(M))\times\gr_{\fv}(N), X\mapsto
(X\cap M, X+N/N),$$
is a surjective morphism of algebraic varieties with affine fibers.
\end{proof}

Now, we can  calculate these characteristics for some particular $Q$-representations.
\begin{equation}\label{simple}
\chi_c(\gr_{\ev}(S_i))=\begin{cases} 1&\hbox{ if } \ev=\dimv(S_i)\hbox{ or } 0\\0&\hbox{ otherwise }\cr\end{cases}  
\end{equation}
Let $P_i$ be the indecomposable projective module with top $S_i$. It is known that if $M$ is a submodule of $P_i$, then either $M=P_i$ or $M\subset$Rad$P_i=\oplus_{i\rightarrow j}P_j$. By this argument and the proposition above, we obtain

\begin{equation}\label{projective}
\chi_c(\gr_{\ev}(P_i))=1+\sum_{\sum\ev^j=\ev}\prod_{ j\rightarrow i}\chi_c(\gr_{\ev^j}(P_j)).
\end{equation}

This together with equation \ref{simple} provides an inductive formula for the ordering of $Q_0$, which enables to compute  the characteritics of projective modules.
For injective modules, which are dual to projective modules, it is sufficient to see that
$$\gr_{\ev}^QN\simeq\gr_{\mv-\ev}^{Q^{opp}}DN,$$
where $D$ is the duality functor from Rep$kQ$ to Rep$kQ^{opp}$.
\end{subsection}

\begin{subsection}{}
We now provide an inductive formula for the Euler characteristics of quiver Grassmannians of a wider class of modules called preprojective and postinjective modules. For this, we need to recall the definition of an almost split sequence and the Auslander-Reiten theorem, \cite{AR}. 

A short exact sequence
\begin{equation}\label{almost split}
0\rightarrow M\rightarrow E\stackrel{\sigma}{\rightarrow} N\rightarrow  0
\end{equation}
in Rep$kQ$ is called {\it almost split} if the following conditions are satisfied
\begin{itemize}
  \item[(i)] $M$ and $N$ are indecomposable
  \item[(ii)] The map $\sigma$ does not split
  \item[(iii)] Given any representation $N'$ and any morphism $\rho$ : $N'\rightarrow N$ which is not a split 
  epimorphism, then $\rho$ factors through $\sigma$.
\end{itemize}  

The proposition that follows is a slight variation of proposition \ref{split}. Actually, the hypo\-the\-sis passes from "split" to "almost split".
\begin{proposition}\label{almost split prop}
Let $0\rightarrow M\rightarrow E\rightarrow N\rightarrow 0$ be an almost split sequence and let $\gv$ be a dimension vector. Then,
$$\chi_c(\gr_{\gv}(E))=\begin{cases}\sum_{\ev+\fv=\gv}\chi_c(\gr_{\ev}(M))\chi_c(\gr_{\fv}(N))&\hbox{ if } \gv\not=\dimv N\\\sum_{\ev+\fv=\gv}\chi_c(\gr_{\ev}(M))\chi_c(\gr_{\fv}(N))+1&\hbox{ if } \gv=\dimv N\cr\end{cases}$$
\end{proposition}
\begin{proof}
As in the proof of proposition \ref{split}, it is sufficient to remark, see \cite[Lemma 3.11]{caldchap}, that the almost split property implies that the map
$$\pi\,:\,\gr_{\gv}(E)\rightarrow\coprod \gr_{\ev}(M))\times\gr_{\fv}(N), X\mapsto
(X\cap M, X+M/M),$$
is a  morphism of algebraic varieties such that \par\noindent
1. case $\gv=\dimv N$: $\pi^{-1}(M',N')$ is empty if $(M',N')=(0,N)$, and is an affine space if not. \par\noindent
2. case $\gv\not=\dimv N$: $\pi^{-1}(M',N')$ is an affine space for any $(M',N')$. 
\end{proof}

The Auslander-Reiten theorem asserts that for any non injective indecomposable representation $M$, there exists a unique almost split sequence beginning with $M$. Dually, for any non projective indecomposable representation $N$, there exists a unique almost split sequence terminating in $N$.

For any non injective indecomposable representation $M$, let $\tau_-M$ be the unique indecomposable  representation (up to isomorphism) $N$ such that \ref{almost split} is almost split.  By convention set $\tau_-I=0$ if $I$ is injective. The indecomposable representations of the form $\tau_-^m(P_i)$, $m\in\N$, $i\in Q_0$ are called (indecomposable) preprojective representations.

Dually, for any non projective indecomposable representation $N$, let $\tau_+N$ be the (unique) indecomposable representation $M$ such that \ref{almost split} is almost split.  By convention set $\tau_+P=0$ if $P$ is projective. The indecomposable representations of the form $\tau_+^m(I_i)$, $m\in\N$, $i\in Q_0$ are called (indecomposable) postinjective representations.

Let ${\mathcal P}$, resp. ${\mathcal I}$, be the set of indecomposable preprojective, resp. postinjective, representations. By the Auslander-Reiten algorithm, ${\mathcal P}$ can be endowed  with an ordering $\preceq$ such that:  
\begin{itemize}
\item[(i)] $\{P_i,\,i\in Q_0\}$ is the set of minimal elements for $\preceq$,\par\noindent
\item[(ii)]  $\preceq$ is generated by $M\preceq\tau_-M$, $M\in{\mathcal P}$
\item[(iii)]  For any $M$ in ${\mathcal P}$, and any summand $E_j$ of $E$ in the almost split sequence $0\rightarrow M\rightarrow E\rightarrow \tau_-(M)\rightarrow 0$, we have $M\preceq E_j\preceq \tau_-M$.
\end{itemize}
Hence, Proposition \ref{almost split prop} gives a recursive formula to calculate the Euler characteristics of any quiver Grassmannian
of a preprojective or postinjective representation. Namely, with the notation above, $\chi_c(\gr_{\fv}(\tau_-M))=$
\begin{equation}
\begin{cases}\chi_c(\gr_{\fv}(E))-\sum_{\ev+\fv'=\fv}\chi_c(\gr_{\ev}(M))\chi_c(\gr_{\fv'}(\tau_-M))&\hbox{ if } \fv\not=\dimv\tau_-(M)\\
1&\hbox{ if } \fv=\dimv\tau_-(M)\end{cases}
\end{equation}
Note that the positivity of Euler characteristics of quiver Grassmannians is not clear from this formula.
It will be the subject of the two next sections.

\end{subsection}
\section{Positivity}\label{positivity}
\begin{subsection}{}
Let $X$ be a complex quasi-projective variety, and let $X_R$ be a model of $X$ over some finitely generated subring $\Z\subset R\subset\C$, that is, $\C\times_R X_R\simeq X$. Given a maximal ideal $m\subset R$, we can form the reduction $X_k$ of $X_R$ to a scheme over the finite field $k=R/m$. We say that $X$ has the polynomiality property if there exists a polynomial $P(t)\in{\bf Z}[t]$ such that, for infinitely many maximal ideals $m$ of $R$, the number of $k$-rational points of $X_k$ equals $P(|k|)$ (note that this property in fact depends on the chosen model $X_R$!). By a straightforward generalization of \cite[Proposition 6.1]{Rrp}, it is sufficient to require existence of a rational function $P(t)\in{\bf Q}(t)$ with this property.\\[1ex]
In particular, we can consider the polynomiality property for the quiver Grassmannian, since it is defined over ${\bf Z}$.
We begin with a proposition which claims that the polynomiality of the quiver Grassmannian implies its "positive polynomiality".
\begin{proposition}\label{hall-positive}
Let $M$ be a finite-dimensional $Q$-representation and let $\ev$ be in $\Z^{Q_0}$. Suppose that the cardinality
of $\gr_{\ev}(M)_{\F_q}$ is a non zero integral polynomial $P_{\ev,M}$ evaluated in $q$. Then, the polynomial $P_{\ev,M}$ has positive coefficients. In particular, we have $\chi_c(\gr_{\ev}(M))>0$.
\end{proposition}
\begin{proof}
The last assertion comes from the classical result $\chi_c(\gr_{\ev(M)})=P_{\ev,M}(1)$, see for example \cite[Lemma 3.5]{caldchap}.

So, we just need to prove that $P_{\ev,M}$ has positive coefficients. Actually, we will show that this is a consequence of Lusztig's theorems on canonical bases of quantum groups, \cite{lusztigbook} and \cite{lusztig}.

For $\mv$ in $\N^{Q_0}$, and $V_{\mv}=\oplus \fq^{m_i}$,
let $\ch_{\mv}(q)$ be the $\C$-vector space of all $G_{V_{\mv}}$-invariant
functions on $E_{\mv}:=E_{V_{\mv}}$.  

Let $\ch(q)$ be the vector space
$\oplus_{\mv\in\N^{Q_0}}\ch_{\mv}(q)$. We can define on $\ch(q)$ a
$\C$-algebra structure by
\[(f_1*f_2)(M)=\sum_{N\subset M}f_1(M/N)f_2(N),\]
where $N$ runs over the (finite) set of submodules of $M$.  This defines an associative unital $\N^{Q_0}$-graded algebra, the Hall algebra of $Q$ (see \cite{Ringel}).

For any isoclass $X$ of $\fq Q$-modules, let $1_X\in\ch(q)$ be the
characteristic function of the corresponding orbit ${\mathcal O}_X\in E_{\dimv X}$. 

We have
\[(1_X*1_Y)(M)=\#\{N\subset M,\,N\simeq Y,\,M/N\simeq X\}.\]
For any dimension vector $\mv$, let $1_{\mv}$ be the constant function with
value 1 on $E_{\mv}$.  Set $\mv= \dimv M$. Then it is easily seen that
\begin{equation}\label{hallpol}
(1_{\mv-\ev}*1_{\ev})(M)=\#\{N,\,N\subset M,\,\dimv N=\ev\}=\#\gr_{\ev}(M)=P_{\ev,M}(q),
\end{equation}

Recall that we fixed an ordering of $Q_0=[1,n]$ such that if
there exists a (non trivial) path from $j$ to $i$, then $j<i$. 

Consider the divided power
\[1_{\alpha_i}^{(e_i)}=\frac{1}{[e_i]_q!}1_{\alpha_i}^{e_i},\]
where $[e_i]_q!=[e_i]_q[e_i-1]_q\ldots [1]_q$, and
$[a]_q:=\sum_{i=0}^{a-1}q^i$ for any integer $a$.  Then we have
\[1_{\alpha_1}^{(e_1)}*1_{\alpha_2}^{(e_2)}\ldots*1_{\alpha_n}^{(e_n)}=1_{\ev}.\]

It is known that $1_{\alpha_i}^{(e_i)}$ is an element of Lusztig's
canonical basis $B$ of $\ch(q)$. Hence, $1_{\ev}$ is a product of elements of the canonical basis, and so is 
$1_{\mv-\ev}*1_{\ev}$.

By Lusztig's positivity theorem, \cite[Theorem 14.4.13]{lusztigbook}, the product of
elements of the canonical basis has positive polynomial coefficients in  the
canonical basis. We can write:
\begin{equation}\label{positivQ}
1_{\mv-\ev}*1_{\ev}=\sum_{b\in B}Q_b(q)b, \,Q_b(q)=\sum_{b,m}c_{b,m}q^m,\,c_{b,m}\in\Z_{\geq 0}.
\end{equation}

Now, in order to apply \ref{hallpol}, we need a description of the evaluation $b(M)$ of the function $b$ on the module $M$. 

A theorem of Lusztig, \cite[Theorem 5.2]{lusztig}, gives such a description. Actually, the function $b$ is obtained by taking an alternating sum of traces of the Frobenius map on the cohomology of a simple perverse sheaf $P_b$ which is pure of weight 0. More precisely, we have

\begin{equation}\label{bm}
b(M)=\sum_i\sum_{j=1}^{\dim H^i(P_b)_M}(-1)^ia_{ij,b,M},\,\hbox{with }\mid a_{ij,b,M}\mid=q^{i/2},
\end{equation}
where $a_{ij,b,M}$, $1\leq j\leq \dim H^i(P_b)_M$, is the set of eigenvalues of the Frobenius map on the stalk at $M$ of the $i$-th cohomology of the $l$-adic complex $P_b$. 

Now, we can use \ref{hallpol}, \ref{positivQ}, \ref{bm} to obtain

\begin{equation}\label{eq1}
P_{\ev,M}(q)=\sum_{b,i,j,m} (-1)^ic_{b,m}q^ma_{ij,b,M}.
\end{equation}

For any nonnegative integer $s$, we can replace $q$ by $q^s$ and obtain

\begin{equation}\label{eqn}
P_{\ev,M}(q^s)=\sum_{b,i,j,m} (-1)^ic_{b,m}q^{ms}a_{ij,b,M}^s.
\end{equation}

As $M$ and $\ev$ are fixed, we now omit $\ev$ and $M$ in the formulas. We set

\begin{equation}
P(q)=P_{\ev,M}(q)=\sum_kp_kq^k.
\end{equation} 

Let $t$ be an indeterminate, we calculate by \ref{eqn} the formal serie $\exp\sum_{n\geq 1}P(q^s)\frac{t^s}{s}$.
As, $p_k\in\Z$ (by hypothesis), $c_{b,m}\in\Z_{\geq 0}$, we obtain:

\begin{equation}
\prod_k\frac{1}{(1-q^kt)^{p_k}}=\prod_{i,j,b,m}\frac{1}{(1-q^ma_{ij,b}t)^{(-1)^ic_{b,m}}}.
\end{equation}

We conclude that  any $a_{ij,b}$ which is not cancelled in the right hand product is an integral power of $q$. 

This implies that any $a_{ij,b}$ which is not cancelled in the right hand sum of \ref{eq1} is an integral power of $q$. In particular, by the last equality of \ref{bm}, we obtain that  if $a_{ij,b}$ is not cancelled in the right hand sum of \ref{eq1}, then $i$ is even.  

Using \ref{eq1} again, we see that $P$ has nonnegative coefficients. 

\end{proof}

We will see in section \ref{reformulation} that this proposition implies the positivity of $\chi(\gr_{\ev}(M))$ in the Dynkin and extended Dynkin case. Recall that the Dynkin case was also proved in \cite{caldkell}.

\end{subsection}
\begin{subsection}{}
We study here the case where the module $M$ has no self-extension. We will prove that the polynomiality, and so the positivity,  is implied by this hypothesis. This example will be fundamental in section \ref{cluster}.

 In the sequel, a $Q$-representation will be called {\it exceptional} if $\Ext^1_Q(M,M)=0$. For example, any indecomposable preprojective or postinjective module is exceptional.
\begin{theorem}\label{except}
Let $Q$ be an acyclic quiver and let $M$ be an exceptional $Q$-representation. Then, for any dimension vector $\ev$ such that $\gr_{\ev}(M)$ is non empty, the cardinality of $\gr_{\ev}(M)\mid_{\F_q}$ is a positive polynomial in $q$. In particular, $\chi_c(\gr_e(M))>0$.

\end{theorem}
\begin{proof}
Set $\mv:=\dimv M$ and $G_{\mv}:=G_{V_{\mv}}$.

From the equations \ref{hallpol} and \ref{positivQ}, we have
\begin{equation}\label{hallpol2}
\#\gr_{\ev}(M)\mid_{\F_q}=1_{\mv-\ev}*1_{\ev}(M)=\sum_{b\in B}Q_b(q)b(M), \,Q_b(q)\in\Z_{\geq 0}[q].
\end{equation}
Now, we need to prove that there exists a unique $b$ such that $b(M)\not=0$ in the sum. 
Indeed, let $b$ be an element of $B$ which is non zero on $M$. 

The corresponding perverse sheaf $ P_b$ is by construction a $G_{\mv}$-invariant perverse sheaf, 
\cite{lusztig}. By construction, see \cite[8.4]{CG}, its support is a closed $G_{\mv}$-subvariety $S_b$ of $E_{\mv}$ and its restriction on a open Zariski dense $O_b$ of $S_b$ is a $G_{\mv}$-invariant local system.

Now, the restriction of $P_b$ on ${\mathcal O}_M$ is a local system. 
As $S_b$ contains $M$, it contains the $G_{\mv}$-orbit ${\mathcal O}_M$ which is an open dense orbit since $M$ is exceptional. So, $S_b=E_{\mv}$.    
As the orbit  ${\mathcal O}_M$ has a connected stabilizer (on $\C$), then \cite[Lemma 8.4.11]{CG} implies that the local system is the trivial sheaf on  ${\mathcal O}_M$. Now, a perverse sheaf whose support is a vector space and which is  trivial  on a Zariski open subset is the constant sheaf.

This implies that $P_b$ is the constant sheaf and so $b(M)=1$.

Equation \ref{hallpol2} gives  
\begin{equation}\label{hallpol3}
\#\gr_{\ev}(M)\mid_{\F_q}=Q_b(q).
\end{equation}
This ends the proof.

\end{proof}

\begin{corollary}
If $M$ is a preprojective or a postinjective module, then any non empty quiver Grassmannian of $M$ has positive Euler characteristics.
\end{corollary}

\end{subsection}

\section{A reformulation of the polynomiality property and its applications}\label{reformulation}
\begin{subsection}{}
Now, let's consider a completed version of the Hall algebra (compare \cite{Rrp}) and perform some computations with certain generating functions in it. By applying an evaluation map, this yields an identity involving cardinalities of quiver Grassmannians over finite fields.

Define $\widehat{\ch}(q):=\prod_{\mv\in\N^{Q_0}}\ch_{\mv}(q)$. The convolution product on $\ch(q)$ of section \ref{positivity} prolongs to a product on $\widehat{\ch}$.
We also consider the skew commutative formal power series ring $\C_q[[\N^{Q_0}]]=\prod_{\mv\in\N^{Q_0}}\C t^{\mv}$ with product $t^{\mv}\cdot t^{\mv'}=q^{-\langle \mv,\mv'\rangle}t^{\mv+\mv'}$. We define the integral as the linear map $\int:\widehat{\ch}(q)\rightarrow\C_q[[\N^{Q_"0}]]$ given by $\int 1_X=\frac{1}{|{\Aut}(X)|}\cdot t^{\dimv X}$ for all $\F_qQ$-representations $X$. The integral is in fact a $\C$-algebra homomorphism by \cite{Rrp}.

We consider the following special elements of $\widehat{\ch}(q)$:
$$e:=\sum_{\mv\in\N^{Q_0}}1_{\mv},\;\;\; e_M:=\sum_{[X]}|{\Hom}(X,M)|1_X,\;\;\;g_M:=\sum_{[X]}|{\Hom}^0(X,M)|1_X,$$
where the last two sums run over all isomorphism classes of finite dimensional representations of $\F_qQ$. We write ${\Hom}(X,M)$ for the space of $\F_qQ$-homomorphism, and ${\Hom}^0(X,M)$ for the subset of monomorphisms.

\begin{lemma} The identity $g_M\cdot e=e_M$ holds in $\widehat{\ch}(q)$, for all representations $M$ of $\F_qQ$.
\end{lemma}

\begin{proof} By the definitions, we have
$$g_M\cdot e=\sum_{[X],[Y]}|{\Hom}^0(X,M)|\cdot 1_X1_Y=\sum_{[Z]}\sum_{Y\subset Z}|{\Hom}^0(Z/Y,M)|\cdot 1_Z.$$
For each representation $Z$, the inner sum on the right hand side is the cardinality of the set of pairs $(Y,g)$ consisting of a subrepresentation $Y$ of $Z$ and an injective map $g:Z/Y\rightarrow M$. But this set is in bijection to ${\Hom}(Z,M)$ by assigning to a map $f:Z\rightarrow M$ the pair $({\ker}(f),\overline{f}:Z/{\ker}(f)\rightarrow M)$. Thus, $$g_M\cdot e=\sum_{[Z]}|{\Hom}(Z,M|\cdot 1_Z=e_X.$$
\end{proof}

We integrate the three elements:
$$\int e=\sum_{[X]}\frac{1}{|{\Aut}(X)|}\cdot t^{\dimv X}=\sum_{\mv}\frac{|E_{\mv}|}{|G_{\mv}|}t^{\mv},$$
since $1/|{\Aut(X)|}$ can be rewritten as $|\co_X|/|G_{\dimv X}|$. Similarly,
$$\int{e_M}=\sum_{\mv}\sum_{a\in\N}\frac{|E_{\mv}^{a,M}|}{|G_{\mv}|}t^{\mv},$$
where $E_{\mv}^{a,M}$ denotes the subset of $E_{\mv}$ consisting of representations $X$ such that the dimension of ${\Hom}(X,M)$ equals $a$. Finally,
$$\int(g_M)=\sum_{[X]}\frac{|{\Hom}^0(X,M)|}{|{\Aut}(X)|}t^{\dimv X}=\sum_{\mv}|{\gr}_{\mv}(M)|t^{\mv}.$$
To see this, we consider the set of pairs $([X],g)$ consisting of an isomorphism class of a representation $X$ of dimension vector $\mv$,
together with an injective map $g:X\rightarrow M$. Sending such a pair to the image of $g$ induces a surjection to ${\gr}_{\mv}(M)$,
whose fibres are precisely the ${\Aut}(X)$-orbits of injective maps.\\[1ex]
Since the integral is an algebra map, the lemma implies the following identity in $\C_q[[\N^{Q_0}]]$:
$$\sum_{\mv'}|{\gr}_{\mv'}(M)|t^{\mv'}\cdot\sum_{\mv''}\frac{|E_{\mv''}|}{|G_{\mv''}|}t^{\mv''}=\sum_{\mv}\sum_{a\in\N}\frac{|E_{\mv}^{a,M}|}{|G_{\mv}|}t^{\mv}.$$
Comparing coefficients, this yields the following recursive formula:

\begin{proposition} The cardinalities of the quiver Grassmannians are given by:
$$|{\gr}_{\mv}(M)|=\sum_{a\in\N}\frac{|E_{\mv}^{a,M}|}{|G_{\mv}|}-\sum_{\mv'<\mv}q^{-\langle \mv',\mv-\mv'\rangle}|{\gr}_{\mv'}(M)|\cdot\frac{|E_{\mv-\mv'}|}{|G_{\mv-\mv'}|}.$$
\end{proposition}

 The above recursive formula now yields:

\begin{corollary}\label{polprop}
 The quiver Grassmannians ${\gr}_{\mv}(M)$ have the polynomiality property if all the varieties $E_{\mv}^{a,M}$ have the polynomiality property.
\end{corollary}

In particular, this yields the polynomiality property of quiver Grassmannians in the case of Dynkin quivers without reference to the existence of Hall polynomials: in this case, each $E_{\mv}^{a,M}$ is a union of finitely many $G_{\mv}$-orbits in $E_{\mv}$, all of which have the polynomiality property by identifying $\co_X$ with the homogeneous space $G_{\mv}/{\Aut}(X)$.

\end{subsection}
\begin{subsection}{}
We can apply corollary \ref{polprop} and proposition \ref{hall-positive} to the extended Dynkin type.

\begin{proposition}\label{affine}
Let $Q$ be a quiver of extended Dynkin type, let $M$ be a $kQ$-module. Then, for any dimension vector $\ev$ the quiver Grassmannian $\gr_{\ev}(M)$ has the polynomiality property. If  $\gr_{\ev}(M)$
is non empty, then its Euler characteristics is positive.

\end{proposition}

\begin{proof} The last assertion is a consequence of the first one by proposition \ref{hall-positive}. To prove the first one, we use proposition \ref{split} to reduce to the case of indecomposable $M$. We will construct a decomposition of $E_{\mv}$ into a disjoint union of strata, such that 
\par\noindent
1. The number of strata is finite,
\par\noindent
2. each stratum has the polynomiality property,
\par\noindent
3. $\dim{\rm Hom}(\_,M)$ is constant along the stratum.

Points 1, 2, 3 imply that each  variety $E_{\mv}^{a,M}$ has the polynomial property; hence Corollary \ref{polprop} provides the first assertion.\\[1ex]
Recall the classification \cite{DlabRingel} of the indecomposable $kQ$-representations into preprojective, postinjective, regular homogeneous and regular non-homogeneous ones, reflecting the classification of the Auslander-Reiten components. All indecomposables of dimension vector a multiple of the minimal imaginary root $\delta$ are represented in these last two classes. Moreover, any regular homogeneous component contains a unique module with given dimension vector $n\delta$, $n>0$.\\[1ex]
As $M$ is indecomposable, it belongs to one of the A.R. components. We denote by $\mathcal{C}_0$ this component if $M$ is regular homogeneous. If not, we can just set $\mathcal{C}_0=\emptyset$ in the sequel.\\[1ex]
Fix isomorphism classes $P$, $I$, $R$, $C$ of a preprojective, resp. postinjective, resp. regular non-homogeneous representation, resp. $\mathcal{C}_0$-component representation, and fix finitely many partitions $\lambda_1,\ldots,\lambda_k$. We define the stratum $S=S_{P,I,R,C,\lambda_1,\ldots,\lambda_k}\subset E_{\mv}$ as consisting of all representations of the form
$$P\oplus I\oplus R\oplus C\oplus \bigoplus_{i=1}^k\bigoplus_{j=1}^\infty U_{i,j}^{(\lambda_i)_j},$$
where the $U_{i,j}$ are regular homogeneous indecomposables of dimension vector $j\cdot\delta$ belonging to $k$  pairwise different homogeneous components $\mathcal{C}_1,\ldots,\mathcal{C}_k$ (different from $\mathcal{C}_0$), that is, $U_{i,j}$ belongs to $\mathcal{C}_i$ for all $i,j$. Since there are at most finitely many isomorphism classes of preprojective or postinjective or regular non-homogeneous indecomposables of any given dimension vector, this defines a finite stratification of $E_{\mv}$. Hence, point 1 is verified.\\[1ex]
Since the number of isomorphism classes of indecomposable representations of any given dimension vector behaves polynomially in the cardinality of the finite base field by \cite{Kac}, the number of isomorphism classes in each stratum $S$ does, too. So, let $P_S(q)$ be the polynomial counting the number of isomorphism classes of indecomposable representations in $S$. 

By standard facts on morphisms between indecomposables over extended Dynkin quivers \cite{DlabRingel} (in particular, the standard directedness properties and orthogonality of regular components), the isomorphism type of the automorphism group of a representation is constant along each stratum. Thus, the cardinality of orbits is constant, say $\#\mathcal{O}_S$, along each stratum. From these facts, we conclude that $\#S=P_S(q)\#\mathcal{O}_S$. Moreover, an orbit has the polynomial property, so point 2 is proved.

Again by \cite{DlabRingel}, $\dim{\rm Hom}(\_,M)$ is constant along each of these refined strata; this gives point 3 and we are done.
\end{proof}

\end{subsection}

\section{The tangent space of the quiver Grassmannian.}{}
This section adopts some methods from \cite{schofield} (see also \cite[section 4.6]{lebruyn}) and is devoted to the calculation of the tangent space of the quiver Grassmannian $\gr_{\ev}(M)$ {\it considered as a scheme}. As a corollary, we obtain the smoothness of the variety $\gr_{\ev}(M)$ when $M$ is an exceptional module. 
\begin{subsection}{}

Let $M=((M_i)_i,(M_\alpha)_\alpha)$ be a representation of $Q$ over $k$ of dimension vector $\mv$, and let $\ev$ be another dimension vector. We still denote by $E_{\ev}$ the variety of representations of $Q$ over $k$ of dimension vector $\ev$, and by $G_{\ev}$ the group acting on it. 

Fix $k$-vector spaces $U_i$ of dimension $e_i$ for all $i\in I$ and define ${\Hom}(\ev,\mv)$ as the space of $I$ graded linear maps $\prod_{i\in I}{\Hom}_k(U_i,M_i)$.

Consider the subvariety ${\Hom}^0(\ev,M)$ of $R_{\ev}\times{\Hom}(\ev,\mv)$ consisting of pairs of tuples $((U_\alpha)_\alpha,(f_i)_i)$ such that all $f_i$ are injective, and $M_\alpha f_i=f_j U_\alpha$ for all arrows $\alpha:i\rightarrow j$. Thus, ${\Hom}^0(\ev,M)$ parametrizes representations of dimension vector $\ev$ together with an injective $kQ$-morphism to $M$. The group $G_{\ev}$ acts  freely on ${\Hom}^0(\ev,M)$  by $(g_i)_i((U_\alpha),(f_i))=((g_jU_\alpha g_i^{-1})_{\alpha:i\rightarrow j},(f_ig_i^{-1}))$.

\begin{lemma} The quiver Grassmannian $\gr_{\ev}(M)$ is isomorphic to the geometric quotient ${\Hom}^0(\ev,M)/G_{\ev}$.
\end{lemma}
\begin{proof} Consider the map $\phi$ : ${\Hom}^0(\ev,M)\rightarrow \gr_{\ev}(M)$, which maps
an element $((U_\alpha),(f_i))$ to $\oplus_if_i(U_i)\subset M$. 

It is well defined since the relations $M_\alpha f_i=f_j U_\alpha$ assert that $\oplus_if_i(U_i)$ is a submodule of $M$. Moreover $\phi$ is clearly surjective.

Suppose that $((U_\alpha),(f_i))$ and $((U_\alpha'),(f_i'))$ have the same image by $\phi$. Hence, for any $i$, $f_i(U_i)=f_i'(U_i)$ and by injectivity of $f_i$, $f_i'$ for any $x$ in $U_i$ there exists a unique $y$ in $U_i$ such that $f_i'(y)=f_i(x)$. The maps $x\mapsto y$ define an element $g_i$ of $\Aut(U_i)$. This gives $(g_i)_i((U_\alpha),(f_i))=((U_\alpha)',(f_i'))$. Hence each fiber of $\phi$ is a (free) $G_{\ev}$-orbit. This implies the lemma.
\end{proof}
Note that the following proposition generalizes a well-known result for classical Grassmannian of vector spaces.
\begin{proposition}
Let $M$ be a $kQ$-module of dimension vector $\mv$, and let $U$ be a submodule of $M$ of dimension vector $\ev$. Then, the tangent space $T_U(\gr_{\ev}(M))$ at $U$ of the quiver Grassmannian
$\gr_{\ev}(M)$ is isomorphic to $\Hom_{kQ}(U,M/U)$.
\end{proposition}
\begin{proof}
To compute the tangent space $T_U(\gr_e(M))$, we fix a point $((U_\alpha),(f_i))$ in ${\Hom}^0(\ev,M)$ in the fiber $\phi^{-1}(U)$, compute the tangent space $T$ of ${\Hom}^0(\ev,M)$ at this point, and factor it by the image of the differential of the action of $G_e$.

To compute $T$, we perform a calculation with dual numbers: since $R_{\ev}\times{\Hom}(\ev,\mv)$ is just an affine space, an element of the tangent space at the point $((U_\alpha),(f_i))$ looks like $((U_\alpha+\epsilon \zeta_\alpha),(f_i+\epsilon A_i))$, with $((\zeta_\alpha),(A_i))\in R_{\ev}\times{\Hom}(\ev,\mv)$. The condition for this to belong to the tangent space $T$ is therefore:
$$M_\alpha(f_i+\epsilon A_i)=(f_j+\epsilon A_j)(U_\alpha+\epsilon \zeta_\alpha)$$
for all $\alpha:i\rightarrow j$, which yields the condition
\begin{equation}\label{eqcondition}
M_\alpha A_i=f_j\zeta_\alpha+A_jU_\alpha,
\end{equation}
 for all $\alpha:i\rightarrow j$. 
 
 The differential of the action of $G_{\ev}$ is computed by applying the definition of the action to a point $((1+\epsilon x_i)_i)$ of $T_1G_{\ev}$:
$$(1+\epsilon x_i)_i((U_\alpha+\epsilon \zeta_\alpha),(f_i+\epsilon A_i))=$$
$$=(((1+\epsilon x_j)(U_\alpha+\epsilon\zeta_\alpha)(1-\epsilon x_i))_{\alpha:i\rightarrow j},((f_i+\epsilon A_i)(1-\epsilon x_i)))=$$
$$=((U_\alpha+\epsilon(\zeta_\alpha+x_jU_\alpha-U_\alpha x_i))_{\alpha:i\rightarrow j},(f_i+\epsilon(A_i-f_ix_i))).$$

By the above, we arrive at the following formula for the tangent space to the quiver Grassmannian:
$$T_{U}({\gr}_{\mv}(M))=\frac{\{((\zeta_\alpha),(A_i))\,\mid\,M_\alpha A_i=f_j\zeta_\alpha+A_jU_\alpha\}}{\{((x_j U_\alpha-U_\alpha x_i)_{\alpha:i\rightarrow j},(-f_ix_i))\,|\, (x_i)_i\in T_1G_e\}}.$$

To understand this condition better, we can assume without loss of generality that $M_\alpha$ is given as a $2\times 2$-block matrix
$$M_\alpha=\left[\begin{array}{rr}U_\alpha&\xi_\alpha\\ 0&N_\alpha\end{array}\right]$$
for all $\alpha$, and that $f_i=\left[{1\atop 0}\right]$ for all $i$ ($1$ representing the identity matrix). We also write $A_i=\left[{B_i\atop C_i}\right]$ for all $i\in I$. Then the condition \ref{eqcondition} reads
$$\left[\begin{array}{rr}U_\alpha&\xi_\alpha\\ 0&N_\alpha\end{array}\right]\left[{B_i\atop C_i}\right]=\left[{1\atop 0}\right]\zeta_\alpha+\left[{B_j\atop C_j}\right]U_\alpha,$$
yielding the two conditions
$$U_\alpha B_i+\xi_\alpha C_i=\zeta_\alpha+B_j U_\alpha$$
and
$$N_\alpha C_i=C_j U_\alpha$$
for all $\alpha:i\rightarrow j$. The subspace to be factored out reads $$\{((x_j U_\alpha-U_\alpha x_i)_{\alpha:i\rightarrow j}),(\left[{-x_i\atop 0}\right]))\}.$$

Hence, choosing $x_i=-B_i$, we can assume each $B_i$ to be zero. We also see that each $\zeta_\alpha$ is uniquely determined by $B_i$ and $C_i$ by the above formula. So the only remaining choices are the $C_i$, subject to the condition
$N_\alpha C_i=C_jU_\alpha$ for all $\alpha:i\rightarrow j$. This condition just means that $C=(C_i)$ defines a $kQ$-morphism from $U$ to $N=(N_\alpha)_\alpha$, and $N$ can be identified with $M/U$. This proves the proposition.
\end{proof}
\end{subsection}
\begin{subsection}{}\label{tangsp}
We give here some consequences of the tangent space formula.

First recall that, since $kQ$ is hereditary, the Euler form $\langle M,M'\rangle:=\dim\Hom_{kQ}(M,M')-\dim\Ext_{kQ}^1(M,M')$ passes to the Grothendieck group. 
\begin{corollary}
Let $M$ be a $kQ$-module of dimension vector $\mv$, and let $U$ be a submodule of $M$ of dimension vector $\ev$. Then, 
$$\langle\ev,\mv-\ev\rangle\leq\dim T_U(\gr_{\ev}(M))\leq \langle\ev,\mv-\ev\rangle+\dim\Ext^1(M,M).$$
In particular, 
$$\langle\ev,\mv-\ev\rangle\leq\dim \gr_{\ev}(M)\leq \langle\ev,\mv-\ev\rangle+\dim\Ext^1(M,M).$$
\end{corollary}

\begin{proof}
Applying ${\Hom}(\_,M)$ to the short exact sequence $0\rightarrow U\rightarrow M\rightarrow M/U\rightarrow 0$,
we get a surjective map $\Ext^1(M,M)\rightarrow \Ext^1(U,M)$, since $kQ$ is hereditary. 

Applying ${\Hom}(U,\_)$ to the same sequence yields a surjection
$$\Ext^1(U,M)\rightarrow \Ext^1(U,M/U).$$
This implies 
$$\dim\Ext^1(U,M/U)\leq\dim\Ext^1(U,M)\leq\dim\Ext^1(M,M).$$
Hence, the equality $\dim{\Hom}_{kQ}(U,M/U)=\langle \ev,\mv-\ev\rangle+\dim\Ext^1(U,M/U)$ gives the proposition.
\end{proof}

\begin{corollary}
Let $M$ be an exceptional module and let $\ev$ be such that $\gr_{\ev}(M)$ is non empty, then $\gr_{\ev}(M)$ is a smooth projective variety of dimension $\langle\ev,\mv-\ev\rangle.$
\end{corollary}
\begin{proof}

Consider the closed subscheme $\gr_{\ev}(\mv)$ of $E_{\mv}\times\prod_i\gr_{e_i}^{m_i}$ consisting of pairs $(M,(U_i)_i)$ of a representation $M$ of dimension vector $\mv$ and a collection $U_i\subset M_i$ of $e_i$-dimensional subspaces such that $$M_h(U_{s(h)})\subset U_{t(h)}$$ for all arrows $h$. This condition defines a subbundle of the trivial vector bundle with fibre $E_{\mv}$ over $\prod_i\gr_{e_i}^{m_i}$, thus $\gr_{\ev}(\mv)$ is smooth and reduced. Its dimension is easily calculated as
$$\sum_{i}e_i(m_i-e_i)+\sum_he_{s(h)}e_{t(h)}+\sum_h(m_{s(h)}-e_{s(h)})m_{t(h)}=\dim E_{\mv}+\langle \ev,\mv-\ev\rangle.$$
The projection $\pi:\gr_{\ev}({\mv})\rightarrow E_{\mv}$ is proper, and its fibre $\pi^{-1}(M)$ is isomorphic to the scheme $\gr_{\ev}(M)$. Now assume that $\gr_{\ev}(M)$ is non-empty, and that $M$ is exceptional. This means that the image of $\pi$ contains the orbit of $M$, which is open. But since $\pi$ is proper, its image is closed. Thus, $\pi$ is surjective. By the theorem on generic smoothness \cite[III.10.7.]{Ha}, this implies that a generic fibre of $\pi$ is smooth, thus reduced, and of dimension
$$\dim \gr_{\ev}(\mv)-\dim E_{\mv}=\langle \ev,\mv-\ev\rangle.$$
Using again the fact that the orbit of $M$ is open in $E_{\mv}$, we conclude that $\gr_{\ev}(M)$ is a smooth variety of dimension $\langle \ev,\mv-\ev\rangle$.
\end{proof}

Note that this corollary does not imply  the last part of theorem \ref{except} since smooth projective varieties may have non positive characteristics.

\end{subsection}

\section{Application to acyclic cluster algebras}\label{cluster}
\begin{subsection}{}
We recall some terminology on cluster algebras.  The reader can find
more precise and complete information in \cite{fominzelevinsky1} and \cite{fominzelevinsky2}.\par

Let $n$ be a positive integer. We fix the {\it ambient field} ${\cF}:=\Q(x_1,\ldots,x_n)$, where the $x_i$'s are indeterminates. Let
${\mathbf x}$ be a free generating set of $\cF$ over $\Q$ and let
$B=(b_{ij})$ be an $n\times n$ antisymmetric matrix with coefficients
in $\Z$. Such a pair $({\mathbf x},B)$ is called {\it a seed}.\par

Let $({\mathbf u}, B)$ be a seed and let $u_j$, $1\leq j\leq n$,
be in ${\mathbf u}$.  We define a new seed as follows. Let $u_j'$ be
the element of $\cF$ defined by the {\it exchange relation}:
\begin{equation}\label{exchangerelation}
u_ju_j'=\prod_{b_{ij}>0}u_i^{b_{ij}}+\prod_{b_{ij}<0}u_i^{-b_{ij}}.
\end{equation}
Set ${\mathbf u'}={\mathbf u}\cup\{u_j'\}\backslash \{u_j\}$.  Let
$B'$ be the $n\times n$ matrix given by
\[
b_{ik}'=\begin{cases}-b_{ik}&\hbox{if } i=j \hbox{ or } k=j\\
b_{ik}+\frac{1}{2}( \,|b_{ij}|\, b_{jk}+b_{ij}\, |b_{jk}|\,) & \hbox{
otherwise.}\cr\end{cases}
\]
By a result of Fomin and Zelevinsky, $({\mathbf u'},B')=\mu_j({\mathbf
u}, B)$ is a seed.  It is called the {\it mutation} of the seed
$({\mathbf u},B)$ in the direction $u_j$ (or $j$).  We consider all
the seeds obtained by iterated mutations from the seed $({\mathbf x},B)$. The free generating sets
occurring in the seeds are called {\it clusters}, and the variables
they contain are called {\it cluster variables}. By definition, the
{\it cluster algebra} $\ca({\mathbf x}, B)$ associated to the seed
$({\mathbf x},B)$ is the $\Z$-subalgebra of $\cF$ generated by the set
of cluster variables. The graph whose vertices are the seeds and whose
edge are the mutations between two seeds is called the {\it mutation
graph} of the cluster algebra.

Let $Q$ be an acyclic finite quiver, we can associate a matrix $B^Q=(b_{ij}^Q)_{i,j\in Q_0}$ defined by
\[b_{ij}^Q=\begin{cases} \#Q_{1,ij}&\hbox{ if } Q_{1,ij} \hbox{ non empty }\\ -\#Q_{1,ji}&\hbox{ if } Q_{1,ji} \hbox{ non empty }\\ 0&\hbox{ otherwise } \end{cases} \]
where $Q_{1,ij}$ is the set of arrows from $i$ to $j$. The cluster
algebra associated to the seed $({\mathbf x},B_Q)$ will be  denoted
by $\ca(Q)$.

\end{subsection}

\begin{subsection}{}

We denote by ${\mathcal E}_Q$ the set of exceptional indecomposable $Q$-representations. According to \cite{caldkell2}, see also \cite{caldchap}, there exists a bijection $\beta$ from the set ${\mathcal E}_Q$ to the set Cl$_Q$ of cluster variables of the cluster algebra $\ca(Q)$.

The bijection is given by 
$$\beta(M)=X_M=\prod_{i\in Q_0}x_i^{-m_i}\sum_{\ev}\chi_c(\gr_{\ev}(M))\prod_{h\in Q_1}x_{s(h)}^{m_{t(h)}-e_{t(h)}}x_{t(h)}^{e_{s(h)}}.$$

Note first that this formula precises, in the acyclic case, the {\it Laurent phenomenon}, see \cite{fominzelevinsky1}, true in general, and which
asserts that the cluster variables are Laurent polynomials with integer
coefficients in the $x_i$, $1\leq i\leq n$. 

This theorem answers a positivity conjecture of Fomin and Zelevinsky in the acyclic case.
It is a direct consequence of Theorem \ref{except}

\begin{theorem}
Let $Q$ an acyclic quiver and $\ca(Q)$ be the associated cluster algebra. Then, the cluster variables are in  $\Z_{\geq 0}[x_1^{\pm 1},\ldots,x_n^{\pm 1}]$
\end{theorem}
We deduce another  proof of the denominator property in the acyclic case, \cite{caldkell2}, \cite{bmrt}, \cite{hubery}.
\begin{corollary}
Let $M$ be an exceptional indecomposable $kQ$-module with dimension vector $\dimv
M=(m_i)$. Then the denominator of $X_M$ as an irreducible fraction of
integral polynomials in the variables $x_i$ is $\prod_ix_i^{m_i}$.
\end{corollary}
\begin{proof}
By the positivity theorem, $X_M$ is a linear combination with positive
integer coefficients of terms $\prod x_l^{n_l}$, $n_l\in\Z$. These
terms are indexed by the set of dimension vectors of submodules $N$ of
$M$, and for each submodule $N$, we have that
$$n_l=-d_l + \sum_{i\to l} e_i + \sum_{l \to j} (d_j -e_j) .$$
Let us put $d=\dim M$ and $e=\dim N$.
Then it is sufficient to prove that
\begin{itemize}
\item[1.] for all $l$, we have $n_l \geq -d_l$ and
\item[2.] for all $l$, there exists a submodule $N$ such that the
equality holds.
\end{itemize}
Since we have $0\leq e_i \leq d_i$ for all $i$, the first assertion is clear.
For the second one, we simply choose $N$ to be the submodule
of $M$ generated by all the spaces $M_j$
such that there exists
an arrow $l \to j$. Then the terms in $\sum_{i\to l} e_i$ vanish
since $Q$ has no oriented cycles and the terms in $\sum_{l \to j} (d_j-e_j)$
vanish by the construction of $N$.
\end{proof}

\end{subsection}
\section{Remarks on constructivity.}
\begin{subsection}{}
We add here some remarks on constructive calculus of the Euler characteristics of the quiver Grassmannian in the acyclic case. Indeed, the use perverse sheaves gives the positivity of $\chi_c(\gr_{\ev}(M))$  but one can't expect  any combinatorial formula from this approach.

It is natural to ask if the quiver Grassmannian $\gr_{\ev}(M)$ has a cellular decomposition, in particular, a decomposition induced by the Schubert cell decomposition of the (subspaces-) Grassmannian of $M$.

In \cite{caldzel}, the authors follow this idea for the Kronecker quiver and obtain closed formulas for the cluster variables and even for a $\Z$-base of the corresponding cluster formula.

Let's sketch a general approach to the problem.  Let $Q$ be any acyclic quiver and set $H:=kQ$. We consider the multiplicative group $H^*$ of invertible elements of $H$. Let $M$ be an $H$-module with dimension vector $\dimv(M)=\mv$ and let $\phi$ be the group homomorphism $H^*\rightarrow \GL(M)$ resulting from the module structure of $M$.

As $Q$ is acyclic, we obtain that $H^*$ is a solvable, so is $\phi(H^*)$, and this implies that there exists a Borel subgroup $B$ of $\GL(M)$ which contains $\phi(H^*)$. 

Let $\underline\gr_{\varepsilon}(M)$ be the (classical) Grassmannian of 
subspaces of $M$ with dimension $\varepsilon$. If $\mid\ev\mid:=\sum_ie_i=\varepsilon$, then $\gr_{\ev}(M)$ is a
closed subvariety of $\underline\gr_{\epsilon}(M)$. The group homomorphism $\phi$ induces an action of $H^*$ on $\underline\gr_{\ev}(M)$.  We claim that
$\coprod_{\mid\ev\mid=\varepsilon}\gr_{\ev}(M)=\underline\gr_{\epsilon}(M)^{H^*}$. 

Indeed, it is clear that any element of $\gr_{\ev}(M)$ is fixed by
$H^*$. Conversely, as $H^*$ is Zariski dense in $H$, any subspace fixed by $H^*$
is stable by $H$ and so, it belongs to $\coprod_{\mid\ev\mid=\varepsilon}\gr_{\ev}(M)$.

Now consider the Schubert cell decomposition relative to $B$:
$$\underline\gr_{\varepsilon}(M)=\coprod_j X_j,$$ where
$X_j$ is an affine space endowed with a left $B$-action.

Hence, there comes:

\begin{proposition}
We have the following decomposition
$$\coprod_{\mid\ev\mid=\varepsilon}\gr_{\ev}(M)=\coprod_j X_j^{H^*},$$
where $X_j$ are affine spaces with a regular $H^*$-action.
\end{proposition}
The problem is that this decomposition has no reason to be a cellular decomposition since the action of $B$ on the affine cell $X_j$ is in general not affine. But the method may afford nice formulas when  the $Q$-representation $M$ is well understood. 
\begin{subsection}{}
Let us describe some results of \cite{caldzel}, obtained by using a Schubert decomposition.

Let $Q$ be the Kronecker quiver with vertices 1 and 2 and two arrows $\alpha$, $\beta$ from 1 to 2. The isoclasses of  indecomposable $Q$-representations where found by Kronecker, \cite{kronecker}, and are described a follows 

Let $M^n$ be the $Q$-representation defined by $\alpha(u_i)=v_i$, $\beta(u_i)=v_{i+1}$, $1\leq i\leq n$, where $(u_1,\ldots, u_n)$, resp. $(v_1,\ldots,v_{n+1})$, is a base of $M_1^n$, resp. $M_2^n$.

Let $DM^n$ be the dual $Q$-representation of $M^n$, ie the spaces are the dual spaces $(M_2^n)^*$, and $(M_1^n)^*$, and the arrows are just the adjoint morphisms.

Let  $M_{reg}(\lambda)^n$, $\lambda\in k$, be the $Q$-representation defined by $\alpha(u_i)=v_i$, $\beta(u_i)=\lambda v_i+v_{i+1}$, $1\leq i\leq n$, where $(u_1,\ldots, u_n)$, resp. $(v_1,\ldots,v_{n})$, is a base of $M_{reg}(\lambda)_1^n$, resp. $M_{reg}(\lambda)_2^n$.

Then, $M^n$, $DM^n$, $n\geq 0$;  $M_{reg}(\lambda)^n$, $n>0$, $\lambda\in k$, is the complete set of isoclasses of indecomposable $Q$-representations.

Using the method above one finds,
\begin{proposition}\cite{caldzel}
$$\chi_c(\gr_{\ev}(M^n))=\begin{cases}\sum_p{p\choose e_2-n-1+p}{e_1-1\choose n-p-e_1}{n+1-e_1\choose p}&\hbox{ if } e_1\not=0\\{n+1\choose e_2}&\hbox{ if } e_1=0\cr\end{cases},$$
$$\chi_c(\gr_{(e_1,e_2)}(DM^n))=\chi_c(\gr_{(e_2,e_1)}(M^n)),$$
$$\chi_c(\gr_{\ev}(M_{reg}(\lambda)^n))=\sum_p{p\choose e_2-n+p}{e_1\choose n-p-e_1}{n-e_1\choose p}.$$

\end{proposition}
 
\end{subsection}

\end{subsection}


\begin{thebibliography}{99}

\bibitem{AR}
M. Auslander, I. Reiten, S. Smalo.
\newblock Representation theory of Artin algebras.
\newblock Cambridge Studies in Advanced Mathematics, 36, Cambridge University Press, Cambridge, 1997.

\bibitem{bmrt}
A. Buan, R. Marsh, I. Reiten, G. Todorov.
\newblock Cluster and seeds in acyclic cluster algebras,
\newblock with an appendix by A. B. Buan, R. J. Marsh, P. Caldero, B. Keller, I. Reiten , and G. Todorov.
\newblock To appear in Proc. Amer. Math. Soc.
\newblock math.RT/0510359.

\bibitem{caldchap}
P.~Caldero, F. Chapoton.
\newblock Cluster algebras as Hall algebras of quiver representations.
\newblock Comment. Math. Helv., 81, (2006), 595-616.
\newblock math.RT/0410184.

\bibitem{caldkell}
P.~Caldero, B. Keller.
\newblock From triangulated categories to cluster algebras.
\newblock to appear in Inv. Math.
\newblock math.RT/0506018.

\bibitem{caldkell2}
P.~Caldero, B. Keller.
\newblock From triangulated categories to cluster algebras II.
\newblock to appear in Ann. Sc. Ec. Norm. Sup.
\newblock math.RT/0510251.

\bibitem{caldzel}
P. Caldero, A. Zelevinsky.
\newblock Laurent expansions in cluster algebras via quiver representations . 
\newblock To appear in  Mosc. Math. J., special issue in honor of Alexander Alexandrovich Kirillov. 

\bibitem{CG}
N. Chriss, V. Ginzburg.
\newblock Representation theory and complex geometry. 
\newblock BirkhŠuser Boston, 1997.

\bibitem{crawley}
W. Crawley-Boevey.
\newblock On homomorphisms from a fixed representation to a general representation of a quiver. 
\newblock Trans. A.M.S., Vol 348, N$^o$5.

\bibitem{DlabRingel}
V. Dlab, C.M. Ringel.
\newblock Indecomposable representations of graphs and algebras.
\newblock Mem. Amer. Math. Soc. 6 (1976), no. 173.


\bibitem{fominzelevinsky1}
 S. Fomin, A. Zelevinsky.
 \newblock Cluster algebras. I. Foundations.
\newblock  J. Amer. Math. Soc.  15  (2002),  no. 2, 497--529.

\bibitem{fominzelevinsky2}
S. Fomin, A. Zelevinsky.
\newblock Cluster algebras. II. Finite type classification.
\newblock Invent. Math.  154  (2003),  no. 1, 63--121.

\bibitem{Green}
J. A. Green.
\newblock Hall algebras, hereditary algebras and quantum groups.
\newblock Invent. Math. 120 (1995), 2, 361-377.

\bibitem{Ha}
R. Hartshorne.
\newblock Algebraic Geometry.
\newblock Graduate Texts in Mathematics, No. 52. Springer, 1977.

\bibitem{hubery}
A. Hubery.
\newblock Hall polynomials for tame hereditary algebras.
\newblock Preprint.

\bibitem{Kac}
V. G. Kac.
\newblock Root systems, representations of quivers and invariant theory.
\newblock Invariant theory (Montecatini, 1982), 74--108, Lecture Notes in Math., 996, Springer, Berlin, 1983.

\bibitem{kronecker}
L. Kronecker,
Algebraische Reduktion der Scharen bilinearer Formen,
\textsl{Sitzungsberichte Akad. Berlin} (1890), 1225–-1237.

\bibitem{lebruyn}
L. Le Bruyn,
\newblock noncommutative geometry@n,
\newblock neverendingbooks, 2005

\bibitem{lusztigbook}
G. Lusztig.
\newblock Introduction to quantum groups.
\newblock Progress in Mathematics 110, Birkh\"auser, Boston, 1993.

\bibitem{lusztig}
G. Lusztig.
\newblock Canonical bases and Hall algebras.
\newblock Representation theories and algebraic geometry (Montreal, PQ, 1997), NATO Adv. Sci. Inst. Ser. C Math. Phys. Sci., vol. 514, Kluwer Acad. Publ., Dordrecht, 1998, pp. 365Ð399. 

\bibitem{Rrp}
M. Reineke.
\newblock Counting rational points of quiver moduli.
\newblock Preprint 2005. To appear in International Mathematical Research Notices.
\newblock math.RT/0505389.

\bibitem{Ringel}
C. M. Ringel.
\newblock Hall algebras and quantum groups.
\newblock Invent. Math. 101 (1990), no. 3, 583--591.

\bibitem{schofield}
A. Schofield.
\newblock General representations of quivers.
\newblock Proc. London Math. Soc. (3) 65 (1992), 1, 46--64.


\end{thebibliography}
\end{document}